\documentclass[11pt,twoside,a4paper]{article}
\usepackage[latin9]{inputenc}
\usepackage{amsmath,amssymb,amsthm,amscd,hyperref}
\usepackage[]{fontenc}
\usepackage{xy}
\usepackage{enumerate}
\usepackage{graphicx}

\newif\ifpdf\ifx\pdfoutput\undefined\pdffalse\else\pdfoutput=1\pdftrue\fi

\xyoption{all} \numberwithin{equation}{section}
\newtheorem{theorem}[equation]{Theorem}
\newtheorem{corollary}[equation]{Corollary}
\newtheorem{proposition}[equation]{Proposition}
\newtheorem{lemma}[equation]{Lemma}

\theoremstyle{definition}
\newtheorem{assumption}[equation]{Assumption}
\newtheorem{definition}[equation]{Definition}
\newtheorem{notation}[equation]{Notation}

\newtheorem{remark}[equation]{Remark}

\newcommand{\nc}{\newcommand}
\nc{\cH}{\mathcal{H}} \nc{\cG}{\mathcal{G}}
\nc{\cC}{\mathcal{C}}
\nc{\cO}{\mathcal{O}}
\nc{\cI}{\mathcal{I}}
\nc{\cB}{\mathcal{B}} \nc{\cY}{\mathcal{Y}} \nc{\cK}{\mathcal{K}}
\nc{\cX}{\mathcal{X}} \nc{\cS}{\mathcal{S}} \nc{\cE}{\mathcal{E}}
\nc{\cF}{\mathcal{F}} \nc{\cZ}{\mathcal{Z}} \nc{\cQ}{\mathcal{Q}}
\nc{\cN}{\mathcal{N}} \nc{\cP}{\mathcal{P}} \nc{\cL}{\mathcal{L}}
\nc{\cM}{\mathcal{M}} \nc{\cR}{\mathcal{R}} \nc{\cT}{\mathcal{T}}
\nc{\cW}{\mathcal{W}} \nc{\cU}{\mathcal{U}} \nc{\cD}{\mathcal{D}}
\nc{\cJ}{\mathcal{J}} \nc{\cV}{\mathcal{V}}
\nc{\fr}{{\rightarrow}}
\nc{\rd}{red.deg}
\newcommand{\Q}{\mathbb{Q}}
\newcommand{\Z}{\mathbb{Z}}
\newcommand{\R}{\mathbb{R}}

\newcommand{\pr}{\mathbb P}

\newcommand{\sym}{\mbox{\upshape{Sym}}}

\newcommand{\rk}{\mbox{\upshape{rank}}}

\title{New slope inequalities for families of complete intersections}
\author{M.A. Barja\footnote{partially supported by the Spanish State Research Agency  AEI / 10.13039/501100011033, through the Severo Ochoa and Mar\'{\i}a de Maeztu Program for Centers and Units of Excellence in R\&D (project CEX2020-001084-M) and the grant PID2019-103849GB-I00. The author is also partially supported by the AGAUR project 2021 SGR 00603 Geometry of Manifolds and Applications, GEOMVAP.} , L. Stoppino\footnote{partially supported MIUR PRIN 2017 ``Moduli spaces and Lie Theory'', by MIUR, Programma Dipartimenti di Eccellenza (2018-2022) - Dipartimento di Matematica ``F. Casorati'', Universit\`a degli Studi di Pavia and by INdAM (GNSAGA)}}

\begin{document}
\maketitle
\begin{abstract}
We prove $f$-positivity of $\cO_X(1)$ for arbitrary dimensional fibrations over curves $f\colon X\to B$  whose general fibre is a complete intersection. In the special case where the family is a global complete intersection, we prove numerical sufficient and necessary conditions for $f$-positivity of powers of $\cO_X(1)$ and for the relative canonical sheaf. From these results we also derive a Chow instability condition for the fibres of relative complete intersections in the projective bundle of a $\mu-$unstable bundle.
\end{abstract}

\pagestyle{myheadings} \markboth{\small{M. A. Barja, L.
Stoppino}}{\small{\textit{New slope inequalities for families of complete intersections}}}

\maketitle

\section{Introduction}\label{intro}

Let $f\colon X\to B$ be a surjective morphism
with connected fibres from an $n$-dimensional smooth projective variety $X $ to a smooth projective curve $B$, with general fibre $F$. 
Let $L$ be a divisor on $X$.
We say that {\em $L$ is $f$-positive} if the following inequality holds:
\begin{equation}\label{f-pos}
e(L):=L^n-n \frac{L_{|F}^{n-1}}{h^0(F, L_{|F})}\deg f_*\cO_X(L)\geq 0.
\end{equation}
The quantity $e(L)$ is a numerical invariant of the data, originally introduced by Cornalba and Harris in \cite{C-H}, and its positivity has a  profound relation with Chow and Hilbert stability of the pair $(F,L_{|F})$ and it is also related with sheaf stability properties of the vector bundle $ f_*\cO_X(L)$: see \cite{BS3} for a detailed discussion; some of these connections will also emerge from the results of the present article.

The case when $L$  is the relative canonical divisor $K_f=K_X-f^{*}K_B$ of a relatively minimal fibration is of particular interest: inequality (\ref{f-pos}) is called in this case {\it the canonical slope inequality}:
\begin{equation}\label{slopeinequality}
K_f^n\geq n\frac{K_{F}^{n-1}}{p_g(F)}\deg f_*\cO_X(K_f).
\end{equation}

The case $n=2$ was the first case studied. Cornalba and Harris  \cite{C-H} (see also \cite{LS}) and Xiao \cite{xiao} proved  with two different methods the canonical slope inequality for relatively minimal fibred surfaces with fibres of genus greater than or equal to $2$. The importance of the slope inequality can hardly be underestimated: in particular it has consequences both in the classification of surfaces (\cite{xiao}, \cite{pardini}) and in the study of the structure of the ample and nef cones of the moduli space of stable curves (\cite{ACG2}, \cite{GKM}).

For the case of an arbitrary $L$ in dimension 2, we proved in \cite{BS3} that $f$-positivity holds in case that $L$ is a nef line bundle such that the restriction $L_{|F}$ is linearly stable \cite[Def.2.16]{mum}, and that this linear stability on the fibres is the assumption needed in both the Cornalba-Harris and the Xiao methods, and also in a third method due to Moriwaki.
It is worth remarking though that -differently from Cornalba-Harris and Moriwaki's ones- Xiao's method works also  with weaker assumptions than linear stability for $L_{|F}$,
 producing  a (weaker) inequality for the ratio between $L^2$  and $\deg f_*\cO_X(L)$.

For the case $n=3$ many inequalities between $K_f^3$ and $\deg f_*\cO_X(K_f)$ are proved in \cite{Barja3folds} and \cite{ohno} and more recently in \cite{zhang}, still a slope inequality or more generally the $f$- positivity of some class of line bundles is in general not known.

For higher values of $n$ only few results are known.
The most important one is in the original paper of Cornalba and Harris \cite{C-H}: an inequality is proved in any dimension under a Hilbert  stability condition on the map induced by the line bundle on the fibres (Theorem \ref{teo-ch}). A slightly improved  result (assuming  the Chow stability of maps induced on the general fibres) has been proved by Bost in arbitrary characteristic \cite{bost}. Unfortunately the Hilbert or Chow stability on the fibres in dimension $\geq 2$ is not known. In particular for the case of $ L=K_f$ and general fibre $F$ of general type, such a condition is known only asymptotically, i.e. for high enough powers of the canonical sheaf.

Very recently the first named author has proved a bound in any dimension for irregular fibrations (see \cite{barja-preprint} and \cite{Barja-BUMI}); similar results have also been proved by Hu and Zhang in \cite{H-Z}.
Even more recently Codogni, Tasin and Viviani in \cite{CTV} applied Xiao's method to families of K-semistale and KSB-semiastable varieties and their moduli.

In \cite{B-S-hyper} we gave a complete treatment of the particular case where $X$ is a relative hypersurface in a projective bundle $\pi\colon \pr_B(\cE)\to B$, where $\cE$ is a vector bundle over $B$, and $f\colon X\to B$ is the morphism induced by $\pi$.
In that case $f$-positivity of $\cO_X(h)$ for any $h\geq 1$ is equivalent to the canonical slope inequality and both are equivalent to a numerical relation between the class of $X$ in the N\'eron-Severi space ${\rm N}^1(\pr)$ and the slope $\mu(\cE):=\deg \cE/\rk \cE$ of $\cE$. From this we could deduce instability and singularity conditions for the fibres and also for the total space $X$.

\bigskip

In this paper we study fibrations whose general fibre is a complete intersection of arbitrary codimension. The main results of the present paper are three.

Firstly we prove a very general slope inequality for these fibrations, with mild conditions on $L$ and on the general fibre $F$.
\begin{theorem}\label{intro-cor-ch}(Theorem \ref{cor-ch})
Let $X$ be an $n$-dimensional variety with a surjective morphism with connected fibres $f\colon X\rightarrow B$ over a smooth curve $B$.
Let $L$ be a line bundle over $X$ which is relatively ample with respect to $f$.
Suppose that
 the general fibre $F$ is embedded in $\pr^{h^0(F, L_{|F})-1}=\pr^{r-1}$ by $|L_{|F}|$ as the complete intersection of $r-n$ hypersurfaces $Y_i$ of degree $d_i$, such that for any $i=1,\ldots r-n$
\[lct(\pr^{r-1},Y_i)\geq \frac{r}{d_i},\]
where $lct$ is the log canonical threshold of the couple $(\pr^{r-1}, Y_i)$.

Then $L$ is $f$-positive.
\end{theorem}
\begin{remark}
In particular the result holds if the general fibre is the (non necessarily smooth) complete intersection of smooth hypersurfaces.
\end{remark}
This is to our knowledge the most general result holding in any dimension. Its proof is however a simple application of Cornalba-Harris and Bost's results, combined with:
\begin{enumerate}
\item a result of Lee relating the Chow stability of a projective variety with the log canonical threshold of its Chow form \cite{Lee};
\item the simple but remarkable fact that the proper intersection of stable varieties is stable (Proposition \ref{stab-intersection}), that follows from a result of Ferretti \cite{ferretti}.
\end{enumerate}

In the second (and longer) part of the  paper we study thoroughly the case of a codimension $c$ complete intersection $X$ in a relative projective bundle $\pr:=\pr_B(\cE)$ over a curve $B$, where $\cE$ is a rank $r\geq 3$ vector bundle over $B$, together with the morphism  $f\colon X\to B$ induced by $\pi$.
Here not only we study the $f$-positivity of $\cO_X(1)$, but also the positivity of its powers, and of the relative canonical sheaf $\omega_f=\cO_X(K_f)$.
It turns out that in many cases  there is a numerical inequality governing the $f$-positivity of these sheaves: inequality \eqref{eq: principale} below.
The main results are the following:
\begin{theorem}[Theorem \ref{thm: piccolo}, Theorem \ref{thm: hgg balanced}, Proposition \ref{prop: meglio che niente}]\label{thm: main}
Let  $X$ be a complete intersection of $c$ hypersurfaces $X_i \equiv k_iH-y_i\Sigma$,  with $k_i\geq 2$ for $i=1,\ldots ,c$.
Consider the inequality:
\begin{equation}\label{eq: principale}
c\mu(\cE)\geq \sum_{i=1}^c\frac{y_i}{k_i}.
\end{equation}
The following three statements are equivalent:
\begin{itemize}
\item [(A1)] the sheaf $\cO_X(h)$ is $f$-positive for any $h< \min\{k_i\}$;
\item [(A2)] there exists a $h< \min\{k_i\}$ such that $\cO_X(h)$ is $f$-positive;
\item[(A3)] inequality \eqref{eq: principale} holds.
\end{itemize}
Suppose moreover that $X$ is balanced (i.e. $k_i=k$ for any $i=1,\ldots , c$); the following statements hold:
\begin{itemize}
\item[(B1)] if  $\cO_X(h)$ is $f$-positive for $h\gg0$ then \eqref{eq: principale} holds;
\item[(B2)] if \eqref{eq: principale} holds with strict inequality then  $\cO_X(h)$ is strictly  $f$-positive for $h\gg 0$.
\end{itemize}
Suppose that $X$ is balanced and that $r<ck$ (i.e.the fibres of $f$ are of general type). Then we have the following:
\begin{itemize}
\item[(C1)] if  $(c-1)k<r$, then $\omega_f$ is $f$-positive if and only if \eqref{eq: principale} holds;
\item[(C2)] if $k\gg0$ and $c=2,3,4$, then $\omega_f$ is $f$-positive if and only if \eqref{eq: principale} holds;
\item[(C3)] if  $c$ is fixed and $r\gg0$, then $\omega_f$ is $f$-positive  if and only if \eqref{eq: principale} holds.
 \end{itemize}
\end{theorem}
\begin{remark}This is a wide generalization of Theorem 2.8 of  \cite{B-S-hyper}, where we proved that  condition \eqref{eq: principale} is equivalent to the $f$-positivity of $\omega_f$ and of  {\em any} power of $\cO_X(1)$ (see Remark \ref{mah}).
\end{remark}
\begin{remark}In the case $c=r-2$ with $X$ balanced, Enokizono made in \cite{ENO} an explicit computation of the invariants $K_f^2$ and $\deg f_*\omega_f$, proving an equality of the form
$K_f^2=\lambda(r,k)\deg f_* \omega_f$. This equality  in particular implies our last results (C1) in case $c=r-2$, as we discuss in detail in Remark \ref{rem: confronto}.
\end{remark}

\begin{remark}\label{rem-ch}
The first Theorem \ref{intro-cor-ch} above is an extremely general result of $f$-positivity, its assumption being on the general fibres.
However, notice that in the case of a global complete intersection in $\pr$ it does not even imply items (A1) and (A2) of Theorem \ref{thm: main}. Indeed, first of all we do not ask any condition on the singularities of the fibres in Theorem \ref{thm: main}; moreover, in (A1) and (A2) we obtain a sufficient and {\em necessary} condition for $f$-positivity. This, combined with Theorem \ref{intro-cor-ch}, can be used to obtain the following strong theorem in the case of relative complete intersections
\end{remark}
\begin{theorem}[Theorem \ref{ thm: finale}]
Let  $X$ be a complete intersection of $c$ hypersurfaces $X_i \equiv k_iH-y_i\Sigma$,  with $k_i\geq 2$ for $i=1,\ldots ,c$. Suppose that for the general fibre $\Sigma\cong \pr^{r-1}$ we have that 
\begin{equation}\label{eq: lct}lct(\Sigma,X_i\cdot \Sigma)\geq \frac{r}{k_i}\qquad \mbox{for any }\quad i=1,\ldots,c.\end{equation}
Then the following statements hold:
\begin{enumerate}
\item the sheaf $\cO_X(h)$ is   $f$-positive for any $h<\min_i\{k_i\}$;
\item inequality \eqref{eq: principale}  holds.
\end{enumerate} 
\end{theorem}

\medskip

In the third and last part of the paper, following the spirit of our work on hypersurfaces \cite{B-S-hyper}, we use the results obtained to study the cones in the N\'eron-Severi space of cycles ${\rm N}^c(\pr)$. Indeed, equation \eqref{eq: principale} tells us something about the position of the class of $X$ inside ${\rm N}^c(\pr)$: we define a cone $\mathbb B$ in the 2-dimensional space ${\rm N}^c(\pr)$ as follows  (Definition \ref{def: cono}).
\[\mathbb B:=\R^+[H^{c-1}\Sigma]\oplus\R^+[H^{c}-c\mu(\cE)  H^{c-1}\Sigma],\] where $\mu(\cE)=\deg\cE/\rk \cE$ is the slope of the vector bundle $\cE$.

By reformulating in a suitable way a result of Fulger \cite{fulger} -using the so-called virtual slopes of the vector bundle $\cE$- we see that this cone $\mathbb B$ is always intermediate between the Pseudoeffective and the Nef cones in ${\rm N}^c(\pr)$. Then we can reformulate the results obtained, for instance by saying that the class of $X$ lies in the interior of $\mathbb B$ if and only if $\cO_X(1)$ is $f$-positive (Proposition \ref{cone}).

Eventually, using the same reasoning as in \cite{B-S-hyper}, we combine the result of Cornalba-Harris and Bost with Theorem \ref{thm: main} and obtain an instability condition for the fibres of relative complete intersections: 
\begin{theorem}[Corollary  \ref{instability}] 
Let $X\subset \pr$ a relative complete intersection in the projective bundle $\pr=\pr_B(\cE)$ satisfying Assumptions \ref{ass: proper}.
If $\sum_{i=1}^c\frac{y_i}{k_i}>c\mu$ (equivalently $[X]\not\in \mathbb B$), then:
\begin{itemize}
\item[(i)] the fibres of $f$ are Chow unstable with the restriction of $\cO_{\mathbb P ^{r-1}}(h)$ for any $h< \min\{k_i\}$.
\item[(ii)] Assume moreover that $X$ is balanced. Then the fibres of $f$ are Chow unstable with the restriction of $\cO_{\mathbb P ^{r-1}}(h)$ for  $h\gg0$;
\item[(iii)] Assume moreover that $X$ is balanced,  $r<kc$ and (1), (2) or (3) in Proposition \ref{prop: meglio che niente} holds. Then the fibres of $f$ are unstable with respect to their dualizing sheaf.
\end{itemize}
\end{theorem}
This leads to an example of unstable complete intersections of general type with only one point as singularity (Proposition \ref{esempio}). Note that such varieties need to be singular (Remark \ref{rem-referee}), so these examples have the smaller possible singularity (set-theoretically).


\section{$f$-positivity of families of complete intersections}\label{pos-stab}
We work over the complex field $\mathbb C$.
In this section we derive some results on the $f$-positivity and the Chow stability of fibres of a fibration whose general fibres are  complete intersections.
For the definition of Chow stability of a projective variety see \cite{ferretti} and the references therein. From now on, anytime we say (semi)stable we mean Chow (semi)stable.
The main result relating this conditions is due to Bost \cite[Theorem 3.3]{bost} and Cornalba-Harris \cite{C-H} \cite{LS} (see \cite{BS3} for references to similar results).

\begin{theorem}[Bost, Cornalba-Harris]\label{teo-ch}
Let $X$ be an $n$-dimensional variety with a surjective morphism $f\colon X\rightarrow B$ over a smooth curve $B$.
Let $L$ be a divisor over $X$ such that:
\begin{itemize}
\item[(i)] for a general fibre $F$, $L_{|F}$ is very ample;
\item[(ii)] $L$  is relatively nef with respect to $f$.
\end{itemize}
If the general fibre of $f$ is Chow semistable with respect to the immersion induced by  $L_{|F}$ then $L$ is $f$-positive.
\end{theorem}

We now ask ourselves: what can we say about the stability of fibres which are complete intersections?
We now state a very natural stability result, which derives from a formula of R. G. Ferretti \cite[Theorem 1.5]{ferretti}.
This application of Ferretti's result was suggested to the second author by Yongnam Lee.

Let $Y$ and $Z$ be two irreducible subvarieties of $\pr^n=\pr(V^{\vee})$ whose intersection is proper.
Let $Y\cdot Z$ be the intersection cycle of $Y$ and $Z$.
\begin{proposition}\label{stab-intersection}
If $Y$ and $Z$ are semistable then $Y\cdot Z$ is semistable.
If, moreover, at least one among $Y$ and $Z$ is stable then the intersection $Y\cdot Z$ is stable.
\end{proposition}
\begin{proof}
We use the Hilbert-Mumford criterion for stability.
Let us consider a 1-parameter subgroup of $GL(V)$ and let $F$ be the associated weighted filtration  of $V$, with weights $r_i$.
Then, for any subvariety $X\subset \pr^n$, there is a well defined integer $e_F(X)$, which in the notation of \cite{ferretti} is called {\em degree of contact}.
The Hilbert-Mumford criterion says that an irreducible subvariety $T\subset \pr(V^{\vee})$ is semistable (resp. stable) if and only if for any weighted filtration $F$ of $V$ the following inequality holds:
\begin{equation}
\frac{e_F(T)}{(\dim T+1)(\deg T)}\leq \frac{1}{n+1}\sum_{i=0}^{n}r_i \quad\quad (\mbox{resp.}<).
\end{equation}
Choosing $Y$ and $Z$ properly intersecting in $ \pr(V^{\vee})$, Ferretti proves in  \cite[Theorem 1.5]{ferretti} the following ``B\'ezout type''
formula for the degree of contact of the cycle intersection $Y\cdot Z$:
\begin{equation}\label{eq-ferretti}
e_F(Y\cdot Z)= \deg(Y)e_F(Z)+\deg(Z)e_F(Y)-\deg(Y)\deg(Z)\sum_{i=0}^{n}r_i.
\end{equation}
Let us now suppose that $Y$ and $Z$ are semistable.
From the Hilbert-Mumford criterion we have that
\[
\frac{e_F(Y)}{(\dim Y+1)(\deg Y)}\leq \frac{1}{n+1}\sum_{i=0}^{n}r_i, \quad\quad \quad  \frac{e_F(Z)}{(\dim Z+1)(\deg Z)}\leq \frac{1}{n+1}\sum_{i=0}^{n}r_i.
\]
Call $y$ and $z$ the dimensions of $Y$ and $Z$ respectively.
By the properness assumption we have that \[\dim(Y\cdot Z)= y+z-n.\]
We thus  have the following chain of inequalities:
\begin{equation*}
\begin{split}
\frac{e_F(Y\cdot Z)}{(\dim (Y\cdot Z)+1)(\deg (Y\cdot Z))} &= \frac{e_F(Y)}{(y+z-n+1)\deg(Y)}+\frac{e_F(Z)}{(y+z-n+1)\deg(Z)}- \frac{\sum_{i=0}^{n}r_i}{y+z-n+1} \\
&\leq\frac{(z+1)+(y+1)}{(y+z-n+1)(n+1)}\sum_{i=0}^{n}r_i- \frac{\sum_{i=0}^{n}r_i}{y+z-n+1}\\
&= \frac{1}{n+1}\sum_{i=0}^{n}r_i,
\end{split}
\end{equation*}
as wanted. The first equality is Ferretti's formula (\ref{eq-ferretti}), while the  inequality derives from the condition of semistability of $Y$ and $Z$.
The result with strict stability follows by substituting strict inequality for (at least) one of the varieties.
 \end{proof}

Let us now recall the definition of log canonical threshold (see \cite{KM} for reference).
Let $(Y, \Delta) $ be a pair, with $Y$ a normal $\Q$-Gorenstein variety and $\Delta$  a $\Q$-Cartier, $\Q$-divisor on $Y$.
Given any birational morphism $\varphi\colon T\longrightarrow Y$ with $T$ normal, we have
\[
K_T + \varphi^{-1}_*\Delta \equiv \varphi^*(K_Y+\Delta) + \sum a( E_i, Y, \Delta) E_i,
\]
where $ \varphi^{-1}_*\Delta $ is the strict transform of $\Delta$ and the $E_i$'s are the exceptional irreducible divisors associated to $\varphi$.
Then we define the {\em discrepancy} of the couple  $discrep(Y, \Delta)$ to be the infimum of the $a(E, Y, \Delta)$, taken for any birational morphism $\varphi$ and any exceptional irreducible divisor.
The couple $(Y, \Delta)$ is said to be {\em log canonical} (l.c.) if $discrep(Y, \Delta)\geq -1$.
The log canonical threshold of $(Y, \Delta)$ is
\[
lct (Y, \Delta):= \sup \{t>0\, | \,\,(Y, t\Delta) \mbox{ is log canonical}\}\in \Q\cap (0, 1].
\]

The log canonical threshold  is a measure of the singularities of the couple $(Y, \Delta)$: for instance if $Y$ is smooth and  $\Delta$ is reduced and normal crossing, then $lct(Y,\Delta)=1$.

In \cite{Lee}, Lee proved the following beautiful condition for a variety to be Chow semistable in terms of the log canonical threshold of its Chow form.
\begin{theorem}[Lee]\label{thmLEE}
Let $Y$ be an $s$-dimensional variety together with a non-degenerate degree $d$ immersion in $\pr^n$.
Let $Z(Y)\subset G$ be the corresponding Chow variety in the Grassmanian $G:= \mbox{Gr}(n-s-1, \pr^n)$.
Suppose that  the following inequality holds
\begin{equation}\label{cond-lct}
lct(G, Z(Y)) \geq \frac{n+1}{d} \quad \mbox{(resp. }> ).
\end{equation}
Then $Y\subset \pr^n$ is Chow semistable (resp.  Chow stable).
\end{theorem}

We are now ready to state the main theorem of this section.

\begin{theorem}\label{cor-ch}
Let $X$ be an $n$-dimensional variety with a surjective morphism $f\colon X\rightarrow B$ over a smooth curve $B$.
Let $L$ be a line bundle over $X$ which is relatively nef with respect to $f$.
Suppose that for the general fibre $F$ the line bundle $L_{|F}$ is very ample and call $r:=h^0(F, L_{|F})$.
Suppose moreover that one of the following conditions holds:
\begin{enumerate}
\item $F$ is embedded in $\pr^{r-1}$ by $|L_{|F}|$ as the complete intersection of $r-n$ hypersurfaces $Y_i$ of degree $d_i$, such that for any $i=1,\ldots r-n$
\[lct(\pr^{r-1},Y_i)\geq \frac{r}{d_i};\]
\item $F$ is embedded in $\pr^{r-1}$ by $|L_{|F}|$ as a degree $d$ variety such that, with notations as in Theorem \ref{thmLEE},
\[lct(G, Z(Y)) \geq \frac{n+1}{d}.\]
\end{enumerate}
Then $L$ is $f$-positive.
\end{theorem}
\begin{proof}
If assumption (1) holds, then by Lee's result above the embeddings  $Y_i\hookrightarrow \pr^{r-1}$ are Chow semistable for any $i=1,\ldots, r-n$, hence  by Proposition \ref{stab-intersection} the complete intersection of the $Y_i$'s is Chow semistable. We thus can apply Theorem \ref{teo-ch} and conclude the proof. For the case of assumption (2) Theorem \ref{thmLEE} directly proves Chow semistability of the intersection, and the rest follows as above.
\end{proof}

\section{Relative complete intersections}\label{sec-inequalities}
\subsection{Set-up and preliminaries}

In the rest of the paper, we make the following assumptions.
Let $\cE$ be a vector bundle of rank $r\geq 3$ and degree $d=\deg (\det \cE)$  on a smooth projective curve $B$ of genus $b$.

Let $\pr:=\pr_B(\cE)$ be the relative projective bundle of duals, following Grothiendieck's notation.
Let $\pi\colon \pr_B(\cE)\rightarrow B$ be the natural projection. Call $\Sigma$ a general fibre of $\pi$ (which is a $\pr^{r-1}$).

Let $\cO_{\pr}(1)$ be the tautological sheaf on $\pr$, and let $H$ be an associated divisor, so that $\cO_{\pr}(1)\cong\cO_{\pr}(H)$.
Let $c$ be an integer between $1$ and $r-2$. Let $X_1, \ldots X_c \subset \pr$ be relative hypersurfaces of degree $k_i\geq 2$.
As $\mathrm{Pic}(\mathbb P)=\Z[\cO_{\pr}(1)]\oplus\pi^* \mathrm{Pic}(B)$, each $X_i$ is an effective divisor in a linear system of the form $|k_iH-\pi^*M_i|$, where  $M_i$ is a divisor on $B$, say of degree $y_i\in \Z$.

Let $X$ be the scheme theoretic intersection of the $X_i$'s, and $f\colon X\rightarrow B $ the induced fibration.
\begin{assumption}\label{ass: proper}
We  assume that  the intersection $X$ is irreducible and proper, i.e. of dimension $r-c$.
\end{assumption}

\begin{definition}
We shall call $X$  {\em balanced} in case that the $k_i$'s are all equal $k_i=k$ for any $i\in \{1,\ldots c\}$.
\end{definition}
Now we would like to use a celebrated result of Miyaoka-Nakayama \cite{nak} about the positive cones of divisors of $\pr_B(\cE)$. 

Before stating the result let us recall the notion of the Harder-Narasimhan sequence of  a vector bundle $\cE$ (see \cite{H-N}). It is the unique filtration of subbundles
$$0=\cE_0\subset \cE_1\subset \ldots \subset \cE_l=\cE$$
satisfying the following assumptions:
\begin{itemize}
\item for any $i=0, \ldots l$ the sheaf $\cE_i/\cE_{i-1}$ is $\mu$-semistable;
\item if we set $\mu_i:= \mu(\cE_i/\cE_{i-1})$, we have that $\mu_i>\mu_{i-1}$.
\end{itemize}
Note that $\mu_1>\mu(\cE)>\mu_l$, unless $\cE$ is $\mu$-semistable, in which case $1=l$ and $\mu_1=\mu(\cE)$. 
Note moreover that, setting $r_0:=0$, we can express the degree of $\cE$ as a combinations of the $\mu_i$'s and $r_i$'s: 
\begin{equation}\label{eq: grado}
\deg(\cE)=\sum_{i=1}^{l}\mu_i(r_{i}-r_{i-1}).
\end{equation}
\begin{theorem}\label{thm: M-N}[Miyaoka-Nakayama]
Given $k\in \Z^{>0}$, $M$ a divisor on $B$,
\begin{enumerate}
\item The  divisor  $kH-\pi^*M$ is pseudo-effective (i.e. it is a limit of effective divisors) if and only if $(\deg M)/k\leq \mu_1(\mathcal E)$, where $\mu_1(\mathcal E)$ is the first slope of the Harder-Narasimhan sequence of $\mathcal E$.
\item The  divisor  $kH-\pi^*M$ is nef if and only if $(\deg M)/k\leq \mu_\ell (\mathcal E)$, where $\mu_\ell (\mathcal E)$ is the last slope of the Harder-Narasimhan sequence of $\mathcal E$.
\end{enumerate}
\end{theorem}
\begin{remark}
Recall that a vector bundle $\cE$ is called nef if the corresponding tautological sheaf $\cO_{\pr(\cE)}(1)$ is nef. From Theorem \ref{thm: M-N} we see that $\cE$ is nef if and only if the smallest slope $\mu_l$ is greater or equal to $0$. In this case, setting $\mu_{l+1}:=0$, we can reformulate equation \eqref{eq: grado} as follows:
\[
\deg(\cE)=\sum_{i=1}^{l}r_i(\mu_{i}-\mu_{i+1}).
\]
\end{remark}
\begin{remark}
Theorem \ref{thm: M-N} tells us that under our Assumptions \ref{ass: proper}  for any $i=1,\ldots , c$, being $X_i\in |k_iH-\pi^*M_i| $ effective, it necessarily holds the inequality $y_i/k_i\leq \mu_1$.
\end{remark}
\begin{notation}\label{notation-multi}
For a multi-index $I  = \{i_1,i_2,...,i_l\} \subseteq  \{1,2,\ldots,c\}$ of length $|I|=l$, we shall indicate by $k_I$ and $y_I$ the corresponding sums:
\[k_I=\sum_{k=1}^{l} k_{i_k},\quad y_I=\sum_{k=1}^{l} y_{i_k}.\]
For $I=\emptyset$, we set $k_I=y_I=0$.

Moreover, we  call $J=\{1,\ldots, c\}$ the full set of indices, so that $k_J=\sum_{i=1}^ck_i$ and  $y_J=\sum_{i=1}^cy_i$.
\end{notation}

Let $H_X$ be the class of $\cO_X(1)$, let $F=\Sigma\cap X$ be the class of a general fibre of $f$, and $H_F$ the class of $\cO_F(1)$.
Consider  the sheaf $\cO_{X}(1)=j^*\cO_\pr(1)$ on $X$, where $j$ is the natural inclusion $j\colon X\hookrightarrow \pr$.
We want to study the  $f$-positivity of $\cO_X(1)$ and of its powers.

Let us start by writing down the basic inequality that embodies the $f$-positivity of $\cO_X(h)$:
\begin{equation}\label{eq: fposHl}
(hH_X)^{r-c}{h^0(F, hH_{|F})} - (r-c) (hH_{|F})^{r-c-1}\deg f_*\cO_X(h) \geq 0.
\end{equation}
We will use the following well known result, whose proof we give for the reader's convenience.

\begin{proposition}\label{prop: sequence}
Let $g\colon Y\to X$ be a morphism between varieties.
Let
\begin{equation}\label{eq: sequenza}
0\longrightarrow \cF_1 \longrightarrow \cF_2\longrightarrow \ldots \longrightarrow \cF_k\longrightarrow 0
\end{equation}
be an exact sequence of ${\cO}_Y$-modules, with $k\geq 3$.

If $R^ig_*{\cF}_j=0$ for any $i>0,j$ with $i+j\leq k-1$, then there is an exact pushforward  sequence:
\[0\longrightarrow g_*{ \cF}_1 \longrightarrow g_* {\cF}_2\longrightarrow \ldots \longrightarrow g_* {\cF}_k\longrightarrow 0.\]
\end{proposition}

\begin{proof}
For $i=1,\ldots ,k-1$, let $\alpha_i:{\cF}_i \longrightarrow  {\cF}_{i+1}$ be the maps in the exact sequence, and for $i=1,\ldots ,k-2$ consider the short exact sequences
\[0\longrightarrow im\alpha_i\longrightarrow {\cF}_{i+1}\longrightarrow im\alpha_{i+1}\longrightarrow 0,\]
where $im\alpha_1={\cF}_1$ and $im\, \alpha_{k-1}={\cF}_k$. Observe that it is enough to prove that $R^1g_* im \alpha_i=0$ for $i=1,\ldots ,k-2$, which is a consequence of the following Claim.

\smallskip

\noindent {\underline {\em Claim.}} For any $i=1,...,k-2$ and for any $j=1,...,k-1-i$ we have $R^jg_*im\,\alpha_i=0.$

\smallskip

\noindent {\underline {\em Proof of claim.}} We proceed by induction on $i$. For $i=1$ we have $im\,\alpha_1=\cF_1$ and so the statement holds by hypothesis. Assume now that the claim holds for $i<k-2$ and consider the short exact sequence

\[
0\longrightarrow im\,\alpha_i\longrightarrow \cF_{i+1}\longrightarrow im\,\alpha_{i+1}\longrightarrow 0.
\]
By induction hypothesis we have that $R^jg_*im\alpha_i=0$ for $j=1,...,k-1-i$ and by the hypothesis in the Proposition we have $R^jg_*{\cF}_{i+1}=0$ for $j=1,...,k-2-i$. So considering the derived long exact sequence we obtain $R^jg_*im \alpha_{i+1}=0$ for $j=1,\ldots ,k-2-i=k-1-(i+1)$.
\end{proof}

Let us now compute the invariants associated to the  $f$-positivity of $\cO_{X}(h)$, following Notation \ref{notation-multi}.
\begin{proposition}\label{koszul}
With the above assumptions, let $h $ be any integer $\geq 1$.
The following formulas hold:
\begin{enumerate}
\item $(H_X)^{r-c}=\left((\prod_{i=1}^ck_i )d-\sum_{i=1}^{c}(\prod_{j\not=i}k_j)y_i\right)$ ;
\item $(H_{F})^{r-c-1}=\prod_{i=1}^{c} k_i;$
\item $ \rk f_*\cO_X(h)=h^0(F, \cO_{F}(h))= \sum_{i=0}^c\left(\sum_{|I|=i}(-1)^i\binom{h-k_I+r-1}{r-1}\right);$\label{eq: rango}
\item\label{grado} $\deg f_*\cO_X(h)= \sum_{i=0}^c\left(\sum_{|I|=i}(-1)^i\binom{h-k_I+r-1}{r-1}\frac{(h-k_I)d+y_Ir}{r}\right).$
\end{enumerate}
\end{proposition}

Note that  here we use the standard convention that considers equal to zero a binomial of the form $\binom{n}{m}$ when $n<m$.
\begin{proof}
The first two formulas are easily computed by standard intersection theory. Indeed we have
\begin{align*}
&H_X^{r-c}=H^{r-c}X=H(k_1H-y_1\Sigma)\ldots(k_cH-y_c\Sigma)=\prod_{i=1}^{c} k_iH^{r}-\sum_{i=1}^{c}(\prod_{j\not=i}k_j)y_iH^{r-1}\Sigma,\\
&H_F^{r-c-1}=H^{r-c-1}F=H^{r-c-1}X\Sigma=H^{r-c-1}(k_1H-y_1\Sigma)\ldots(k_cH-y_c\Sigma)\Sigma=\prod_{i=1}^{c} k_iH^{r-1}\Sigma,
\end{align*}
and now we just need to recall that $H^{r}=\deg\mathcal E=d$ and $H^{r-1}\Sigma=1$.

For the last two formulas, let us consider the Koszul sequence that provides the resolution of $\cO_X(h)$:
(cf. for instance \cite[Chap.III pp 144-145]{E-H}):
\begin{gather*}
0\longrightarrow \cO_{\pr}(-X_1-\ldots -X_c+hH)\longrightarrow \ldots \longrightarrow  \oplus_{|I|=l}\cO_{\pr}(-X_I+hH)\longrightarrow \ldots\\
 \ldots\longrightarrow  \oplus_{i=1}^c\cO_{\pr}(-X_i+hH)\longrightarrow \cO_\pr (h)\longrightarrow \cO_X(h) \longrightarrow 0,
\end{gather*}
where we use the notation $X_I$ to indicate the divisor $\sum_{j\in I} X_{j}$.
For $1\leq i\leq r-2$,  we have the vanishing of the higher direct image sheaves $R^i\pi_*\cO(-X_I+hH)=0$, so we are in the conditions of applying Proposition \ref{prop: sequence}, and obtain by pushforward via $\pi$ the long exact sequence
\begin{gather*}
0\longrightarrow \pi_* \cO_{\pr}(-X_1-\ldots -X_c+hH)\longrightarrow \ldots \longrightarrow  \oplus_{|I|=l}\pi_* \cO_{\pr}(-X_I+hH)\longrightarrow \ldots\\
 \ldots\longrightarrow  \oplus_{i=1}^c\pi_* \cO_{\pr}(-X_i+hH)\longrightarrow \pi_*\cO_\pr (h)\longrightarrow f_*\cO_X(h) \longrightarrow 0.
\end{gather*}
It is then  straightforward to compute  $\rk f_*\cO_X(h)$ obtaining (\ref{eq: rango}).

In order to prove the last formula, observe that for any integer $a\geq 0$ and for any $M $ divisor on $B$ of degree $y$, the projection formula \cite[Chap.~II, Ex.5.1(d)]{Har} says that
\[ \pi_*\cO_\pr(aH-\pi^*M)\cong \pi_*\cO_\pr (a) \otimes \cO_B(-M)\cong \sym^{a}\cE\otimes  \cO_B(-M),\]
 and that
\[\deg \pi_*\cO_\pr(aH-\pi^*M)=\deg  \sym^{a}\cE\otimes \cO_B(-M) = \binom{a+r-1}{r-1}\frac{ad-yr}{r}.\]
Notice moreover that for $a<0$ we have $\pi_*\cO_\pr(a)=0$, and so the formula above still holds.
An easy  computation now leads to formula (\ref{grado}) for  $\deg f_*\cO_X(h)$.
\end{proof}

\subsection{$f$-positivity of $\cO_X(h)$ for small enough $h$}
We first prove a result about the case of  small $h$.
\begin{theorem}\label{thm: piccolo}
Let $f\colon X\to B$ be a morphism as in Assumption \ref{ass: proper}.

The following statements are equivalent:
\begin{enumerate}
\item\label{eq: num} $\sum_{i=1}^c \frac{y_i}{k_i}\leq c\mu({\cal E})$;
\item\label{low f-pos} the line bundle $\cO_{X}(h)$ is $f$-positive for any $h<\min_i\{k_i\}$;
\item\label{uno} there exists an $h<\min_i\{k_i\}$ such that $\cO_{X}(h)$ is $f$-positive.
\end{enumerate}
\end{theorem}
\begin{proof}
Combining the formulas of Proposition \ref{koszul}, we see that  $f$-positivity of $\cO_X(h)$ is equivalent to the following inequality
\begin{equation*}\begin{split}
h\left(\prod_{i=1}^ck_i d- \sum_{i=1}^{c}(\prod_{j\not=i}k_j)y_i\right)\sum_{i=0}^c\left(\sum_{|I|=i}(-1)^i\binom{h-k_I+r-1}{r-1}\right)\geq  \\
(r-c)\prod_{i=1}^{c} k_i \left[  \sum_{i=0}^c\left(\sum_{|I|=i}(-1)^i\binom{h-k_I+r-1}{r-1}\frac{(h-k_I)d+y_Ir}{r}\right)\right].
\end{split}
\end{equation*}
Grouping terms, this inequality becomes
\begin{equation}\label{eq: ine}
\begin{split}
\frac{h}{r}\left[ c\prod_{i=1}^{c} k_i d-\sum_{i=1}^{c}(\prod_{j\not=i}k_j)y_i r \right] \sum_{i=0}^c\left(\sum_{|I|=i}(-1)^i\binom{h-k_I+r-1}{r-1}\right) +\\
+(r-c) \prod_{i=1}^ck_i \left[  \sum_{i=0}^c\left(\sum_{|I|=i}(-1)^i\binom{h-k_I+r-1}{r-1}\frac{k_Id-y_Ir}{r}\right)\right]\geq 0.
\end{split}
\end{equation}
Now observe that for $h< \min_i\{k_i\}$ (hence in particular for $h=1$), the second term in inequality (\ref{eq: ine}) vanishes when $i>0$ because all binomials are zero, and trivially vanishes when $i=0$. Hence in this case the inequality just becomes
\begin{equation}\label{mezza}
\frac{h}{r}\binom{h+r-1}{r-1}\left(c\prod_{i=1}^{c} k_i d-\sum_{i=1}^{c}(\prod_{j\not=i}k_j)y_i r \right)\geq 0.
\end{equation}
From these observations we see that  the condition for $f$-positivity for $h<\min_i\{k_i\}$ is precisely
\begin{equation}\label{eq: condition}
c\prod_{i=1}^{c} k_i d-\sum_{i=1}^{c}(\prod_{j\not=i}k_j)y_i r\geq 0\quad  \Longleftrightarrow \quad \sum_{i=1}^c \frac{y_i}{k_i}\leq c\frac{d}{r}= c\mu,
\end{equation}
thus proving the equivalence between (\ref{eq: num}) and (\ref{low f-pos}). Observe now that if condition (\ref{uno}) holds, then inequality \eqref{mezza} holds, so it holds \eqref{eq: condition}, so (\ref{eq: num}) and (\ref{low f-pos}) both hold.
Note that the same statement holds for {\em strict} $f$-positivity with strict inequality in \eqref{eq: ine}.
\end{proof}
\begin{remark}
This result in particular implies that strict $f$-positivity for relative complete intersections is stable under intersection. More precisely, consider two relative complete intersections $X$ and $X'$ (of codimension $c$ and $c'$ respectively and with associated string of integers $k_i, y_i$ and ${k'}_i, {y'}_i$ respectively, and call $f$ and $f'$ the morphisms induced by $\pi$) in $\pr$. Suppose that $X$ and $X'$ intersect properly, call $X'':=X\cdot X'$ and call $f''\colon X''\to B$ the induced fibration. Suppose that the tautological bundles $\cO_X(1)$ (resp. $\cO_{X'}(1)$) are strictly $f$-positive (resp. $f'$-positive); then by the above theorem we have
\[
\sum_{i=1}^{c}\frac{y_i}{k_i}+\sum_{i=1}^{c}\frac{{y'}_i}{{k'}_i}<c\mu(\mathcal E)+c'\mu(\mathcal E)=(c+c')\mu(\mathcal E),
\]
so, $\cO_{X''}(1)$ is strictly $f''$-positive. This property should be compared with the analogue result holding for Chow stability (Proposition \ref{stab-intersection}).
\end{remark}

\subsection{$f$-positivity of $\cO_X(h)$ for big enough $h$}
\begin{remark}\label{rem: mmm}
Let us consider the following number associated to our data:
\[
\alpha:=c\prod_{i=1}^{c} k_i d-\sum_{i=1}^{c}(\prod_{j\not=i}k_j)y_i r.
\]
In the pervious section we have seen that $f$-positivity of low powers of  $\cO_{X}(h)$ is related to the positivity of $\alpha$. 

Note moreover  that if $\frac{y_i}{k_i}\leq  \mu({\cal E})$ (resp. $<$) for any $i=1,\ldots c$ then $\alpha \geq 0$ (resp $>$). 
\end{remark}
Let us start by recalling the asymptotic formula of Hirzebruch-Riemann-Roch  \cite[Chap. 15 and 18]{fulton}:
\begin{proposition}[Hirzebruch-Riemann-Roch]
Let $L$ be an ample line bundle over a smooth $n$-dimensional variety $Z$. 
We have for $h\gg0$
\begin{equation}\label{eq: asy hrr}
h^0(Z,\cO_Z(hL))=\frac{h^n}{n!}L^n-\frac{h^{n-1}}{2(n-1)!}L^{n-1}K_Z+\cO (h^{n-2}).
\end{equation}
The same formula holds true if $Z$ is a (non necessarily smooth) local complete intersection in a projective space.
\end{proposition}

\begin{theorem}\label{thm: hgg}
Let $f\colon X\to B$ be a morphism as in Assumption \ref{ass: proper}.
The following statements hold:
\begin{itemize}
\item[(1)] Assume that  $\frac{y_i}{k_i}\leq  \mu({\cal E})$ for any $i=1,\ldots c$ and that for some $j\in \{1,\ldots c\}$ strict inequality holds;
then the line bundle $\cO_{X}(h)$ is strictly $f$-positive for $h\gg0$.
\item[(2)] Assume that $\frac{y_i}{k_i}>  \mu({\cal E})$ for any $i=1,\ldots c$; then  for $h\gg0$ the line bundle $\cO_{X}(h)$ is not  $f$-positive.
\end{itemize}
\end{theorem}
\begin{proof}
Let us start from the equation of $f$-positivity (\ref{eq: fposHl}):
\begin{equation}\label{eq: f-pos eq}
\begin{split}
 h^{r-c}H_X^{r-c}\rk f_*\cO_X(h)-(r-c)h^{r-c-1}H_F^{r-c-1}\deg f_*\cO_X(h) \geq 0
\end{split}
\end{equation}
By using Hirzebruch-Riemann-Roch, and the computations in Proposition \ref{koszul}, we have that the left side term is for $h\gg0$ a polynomial in $h$ of the following form:
\begin{equation*}
\begin{split}
 h^{r-c}H_X^{r-c}\rk f_*\cO_X(h)-(r-c)h^{r-c-1}H_F^{r-c-1}\deg f_*\cO_X(h)= \\
 h^{r-c-1}\left(h H_X^{r-c}h^0(F,\cO_F(h))-(r-c)H_F^{r-c-1}\deg f_*\cO_X(h)\right)=\\
 h^{r-c-1}\left(h^{r-c}H_X^{r-c}\frac{H_F^{r-c-1}}{(r-c-1)!}-(r-c) h^{r-c-1}H_F^{r-c-1}\frac{H_X^{r-c}}{(r-c)!}\right)+\cO(h^{2r-2c-2})=\\
 h^{2r-2c-1}\left( \frac{H_X^{r-c}H_F^{r-c-1}}{(r-c-1)!}-(r-c) \frac{H_X^{r-c}H_F^{r-c-1}}{(r-c)!}  \right) +\cO(h^{2r-2c-2}).&
\end{split}
\end{equation*}
So we see that the degree is at most $2r-2c-2$.
We will compute the coefficient of $h^{2r-2c-2}$.

Using  the third and fourth formulas of Proposition \ref{koszul}, grouping terms as in the proof of Theorem \ref{thm: piccolo} above, we see that for any $h\in \mathbb N$ we have:
\[
\deg f_*\cO_X(h)= \frac{dh}{r}\rk f_*\cO_X(h)-\left[  \sum_{i=0}^c\left(\sum_{|I|=i}(-1)^i\binom{h-k_I+r-1}{r-1}\frac{k_Id-y_Ir}{r}\right)\right].
\]
We want to see the right hand term as a polynomial in $h$ and compute its leading term.
\begin{lemma}\label{lem: asi}
We have:
\begin{equation*}
\begin{split}
\beta (h):=&\left[  \sum_{i=0}^c\left(\sum_{|I|=i}(-1)^i\binom{h-k_I+r-1}{r-1}\frac{k_Id-y_Ir}{r}\right)\right]\\
&=-\frac{\alpha}{r}\frac{h^{r-c}}{(r-c)!}+\frac{h^{r-c-1}}{2(r-c-1)!}\left(\frac{\alpha}{r} (k_J-r)+\prod_{i=1}^ck_i\sum_{i=1}^c\frac{k_id-y_ir}{r}\right) +\cO(h^{r-c-2}).
\end{split}
\end{equation*}
\end{lemma}
\begin{proof}{\em of the Lemma.}
Let us compute the coefficient of $\frac{k_1d-y_1r}{r}$ in $\beta(h)$:
it is immediate to see that this is
\[
\sum_{i=1}^c\left(\sum_{1\in I, |I|=i}(-1)^i\binom{h-k_I+r-1}{r-1}\right)=-\sum_{j=0}^{c-1}\left(\sum_{1\not\in J, |J|=j}(-1)^j\binom{h-k_1-k_J+r-1}{r-1}\right),
\]
and this, by the very same computation as in Proposition \ref{koszul}, is precisely \[-\rk {f_1}_*\cO_{Y_1}(h-k_1),\] where $Y_1$ is the $(c-1)$-codimensional intersection of all the $X_i$'s except for $X_1$, and $f_1\colon Y_1\to B$ is the induced fibration. Let us call $H_1$ the tautological divisor on $Y_1$ and $F_1$ the class of a general fibre of $f_1$.
Now, by Hirzebruch-Riemann-Roch \eqref{eq: asy hrr}, this rank has the following asymptotic expansion for $h\gg 0$:
\begin{align*}
 \rk {f_1}_*\cO_{Y_1}(h-k_1)  &=h^0(F_1, (h-k_1)H_{F_1})\\
 &=(h-k_1)^{r-c} \frac{(H_{|F_1})^{r-c}}{(r-c)!} -(h-k_1)^{r-c-1}\frac{(H_{F_1})^{r-c-1}K_{F_1}}{2(r-c-1)!} +\cO(h^{r-c-2})\\
&=h^{r-c} \frac{(\prod_{i\not =1 }k_i)}{(r-c)!} - h^{r-c-1}\left(\frac{k_1(r-c)H_{F_1}^{r-c}}{(r-c)!}+ \frac{H_{F_1}^{r-c-1}K_{F_1}}{2(r-c-1)!}\right)+\cO(h^{r-c-2})\\
&=h^{r-c} \frac{(\prod_{i\not =1 }k_i)}{(r-c)!}- h^{r-c-1}\left(\frac{2k_J+\prod_{i\not=1} k_i(\sum_{i\not=1}k_i-r)}{2(r-c-1)!}\right)+\cO(h^{r-c-2}),
\end{align*}
where we used that:
\begin{itemize}
\item $H_{F_1}^{r-c}=\prod_{j\not=1}k_j$;
\item $H_{F_1}^{r-c-1}K_{F_1}=\prod_{j\not=1}k_j(\sum_{j\not=1}k_j-r),$
\end{itemize}
because $K_{F_1}\equiv (\sum_{j\not=1}k_j-r)H_F$.
So, for $h\gg 0$ we obtain from the above computations a term of the form
\[
-(\prod_{i\not =1 }k_i) \frac{k_1d-y_1r}{r}=-\frac{(\prod_{i=1}^{c} k_i )d-(\prod_{j\not=1}k_j)y_1r }{r}.
\]
Now, recalling that  
\[
\beta(h)=\sum_{i=1}^c\left[  -\left(\frac{k_id-y_ir}{r}\right)\rk f_{1*}\cO_{Y_i}(h-k_i)\right],
\]
we obtain, summing up for any $i=1,\ldots c$, as a degree $r-c$ term for $\beta(h)$
\[
-\frac{h^{r-c}}{(r-c)!}\frac{\alpha}{r}.
\]
Let us compute now the term in degree $r-c-1$. We have 
\begin{equation*}
\begin{split}
\frac{h^{r-c-1}}{2(r-c-1)!}\left[\sum_{i=1}^c\frac{k_id-y_ir}{r}\left( 2\prod_{j=1}^c k_j+ \prod_{j=1, j\not=i}^c k_j(\sum_{j=1, j\not=i}^ck_j-r)\right)\right]\\
=\frac{h^{r-c-1}}{2(r-c-1)!}\left[\sum_{i=1}^c\frac{k_id-y_ir}{r}\prod_{j=1, j\not=i}^c k_j\left( 2k_i+\sum_{j=1, j\not=i}^ck_j-r\right)\right]\\
=\frac{h^{r-c-1}}{2(r-c-1)!}\left[(k_J-r)\frac{\alpha}{r}+(\prod_{j=1}^ck_j)\sum_{i=1}^c \frac{k_id-y_ir}{r}\right],
\end{split}
\end{equation*}
as wanted.
\end{proof}
\noindent Let us now resume the proof of Theorem \ref{thm: hgg}. Consider again the $f$-positivity of $\cO_X(h)$:
\begin{equation*}
\begin{split}
 h^{r-c}H_X^{r-c}\rk f_*\cO_X(h)-(r-c)h^{r-c-1}H_F^{r-c-1}\deg f_*\cO_X(h) \geq 0
\end{split}
\end{equation*}
Dividing by $h^{r-c-1}$ and using  the computations of Proposition \ref{koszul}, it is equivalent to
\begin{equation}\label{eq: uff}
h \left( H_X^{r-c}-\frac{d(r-c)}{r}{H_F}^{r-c-1}\right)\rk f_*\cO_X(h)+(r-c) (\prod_{i=1}^{c} k_i) \beta(h)\geq 0.
\end{equation}
Now observe that by Proposition \ref{koszul} again
\[
H_X^{r-c}-d\frac{r-c}{r}{H_F}^{r-c-1}=\frac{r (\prod_{i=1}^{c} k_i) d- r(\sum_{i=1}^{c}(\prod_{j\not=i}k_j)y_i) -(r-c)d(\prod_{i=1}^{c} k_i)}{r}=\frac{\alpha}{r},
\]
Let us recall Hirzebruch-Riemann-Roch \eqref{eq: asy hrr} for $\rk f_*\cO_X(h)=h^0(F, \cO_F(h))$ for $h\gg 0$:
\begin{equation*}
\begin{split}
\rk f_*\cO_X(h)= h^0(F, \cO_{F}(h)) = \frac{h^{r-c-1}}{(r-c-1)!}H_{F}^{r-c-1} -\frac{h^{r-c-2}}{2(r-c-2)!}H_{F}^{r-c-1}K_F+ \cO(h^{r-c-3})\\
=  \frac{h^{r-c-1}}{(r-c-1)!} (\prod_{i=1}^{c} k_i)-\frac{h^{r-c-2}}{2(r-c-2)!}(\prod_{i=1}^{c} k_i)(k_J-r)+ \cO(h^{r-c-3}),\\
\end{split}
\end{equation*}
because $K_F\equiv (k_J-r)H_F$.
Using the above expression and  and Lemma \ref{lem: asi}, we have that the term in degree $r-c$ in  \eqref{eq: uff} is
\[
\frac{h^{r-c}}{(r-c-1)!}(\prod_{i=1}^{c} k_i)\frac{\alpha}{r}+(r-c)\frac{h^{r-c}}{(r-c)!}(\prod_{i=1}^{c} k_i)\frac{(-\alpha)}{r}= 0,
\]
according to the observation made at the beginning of the proof.

Now,  we compute the degree $r-c-1$ term: 

\begin{equation}\label{eq: termine}
\begin{split}
\frac{h^{r-c-1}}{2(r-c-1)!}(\prod_{i=1}^{c} k_i)\frac{\alpha}{r}(k_J-r)\left( -(r-c-1)+(r-c)\right)+(\prod_{i=1}^{c} k_i)(r-c)\prod_{j=1}^ck_j\sum_{i=1}^c \frac{k_id-y_ir}{r}\\
=\frac{h^{r-c-1}}{2(r-c-1)!}(\prod_{i=1}^{c} k_i)\left((k_J-r)\frac{\alpha}{r}+(r-c)\frac{\gamma}{r}\right),
\end{split}
\end{equation}
where we set 
\[\gamma= \prod_{j=1}^ck_j\sum_{i=1}^c \frac{k_id-y_ir}{r}.\]
So, we see that if we assume $k_J-r>0$ (i.e. $F$ of general type) this coefficient is a strictly positive combination of $\alpha$ and $\gamma$. It is now immediate to check (see also Remark \ref{rem: mmm}) that in the assumption of (1) both $\alpha$ and $\gamma$ are strictly positive. However, we now will see that the assumption $k_J-r>0$  is not needed. Let us rearrange 
the terms in \eqref{eq: termine} as follows:
\begin{equation*}
\begin{split}
\left((k_J-r)\frac{\alpha}{r}+(r-c)\frac{\gamma}{r}\right)=(\prod_{j=1}^{c} k_j)\left((k_J-r)\sum_{i=1}^c \frac{k_id-y_ir}{k_ir}+(r-c)\sum_{i=1}^c \frac{k_id-y_ir}{r}\right)\\
=(\prod_{j=1}^{c} k_j)\sum_{i=1}^c \frac{k_id-y_ir}{k_ir}(k_J-r+ (k_i(r-c))\\
=(\prod_{j=1}^{c} k_j)\sum_{i=1}^c \frac{k_id-y_ir}{k_ir}\left(k_J-c+(k_i-1)(r-c)\right).
\end{split}
\end{equation*}
So, 
it is clear that if the assumption of (1) hold, then the coefficient of $h^{r-c-1}$ is strictly positive, and so $\cO(h)$ is $f$-positive for $h\gg0$.

If on the contrary all the terms $k_id-y_ir$ are strictly smaller than $0$, for $i=1,\ldots c$, then so is  the coefficient of $h^{r-c-1}$ in  \eqref{eq: termine}, and proposition (2) is thus proved.
\end{proof}


We now see that in the balanced case we can trace back the asymptotic  $f$-positivity of $\cO_C(h)$  to the positivity of $\alpha$.
\begin{theorem}\label{thm: hgg balanced}
Let $f\colon X\to B$ be a morphism as in Assumption \ref{ass: proper}. Assume moreover that $f$ is balanced.
The following implications hold.
\begin{itemize}
\item If
 $\frac{y_J}{k}=\sum_{i=1}^c \frac{y_i}{k}< c\mu({\cal E})$, then
 the line bundle $\cO_{X}(h)$ is strictly $f$-positive for $h\gg0$.
\item Conversely, if $\cO_{X}(h)$ is $f$-positive for $h\gg0$ then  $\frac{y_J}{k}=\sum_{i=1}^c \frac{y_i}{k}\leq c\mu({\cal E})$.
\end{itemize}
\end{theorem}
\begin{proof}
The proof of this result is straightforward by computing the coefficient of $h^{r-c-1}$ in  \eqref{eq: termine} in the balanced case:
we obtain 
\[
k^c\sum_{i=1}^c\frac{kd-y_ir}{kr}(ck-r+k(r-c))=(k-1)\alpha.
\]
So, both the implications are clear. 
\end{proof}
\begin{remark}
Clearly  in case of equality $\alpha=0$, in order to check  $f$-positivity, one should investigate the non-negativity of the coefficient in degree ($2r-2c-3$) in the polynomial (\ref{eq: fposHl}).
\end{remark}

\begin{remark}\label{mah}
In case $c=1$  both terms in inequality \eqref{eq: ine} are multiple of the number $(dk_1-ry_1)$ and one can see -as we do in \cite{B-S-hyper}- that the multiplying term is positive, so that in Theorem \ref{thm: hgg} the condition of $\cO_X(h)$ being $f$-positive for $h\gg 0$ is equivalent to $\cO_X(h)$ being $f$-positive for {\em any} fixed value $h\geq 1$.
\end{remark}

Let us now put together the results of Section \ref{pos-stab} with the ones in this section and obtain the following:
\begin{theorem}\label{ thm: finale}
Let $f\colon X\to B$ be a morphism as in Assumptions \ref{ass: proper}.
Suppose that the relative hypersurfaces $X_i$'s are such that  for the general fibre $\Sigma\cong \pr^{r-1}$  of $\pi\colon \pr(\cE)\to B$ we have that
\[lct(\Sigma,X_i\cdot \Sigma)\geq \frac{r}{k_i}.\]
Then 
\begin{enumerate}
\item the sheaf $\cO_X(h)$ is $f$-positive for any $h<\min_i\{k_i\}$;
\item $\sum_{i=1}^c \frac{y_i}{k_i}\leq  c\mu({\cal E})$.
\end{enumerate} 
\end{theorem}
\begin{proof}
The assumption on the log canonical threshold of the couple $(\Sigma, X_i\cdot \Sigma)$ implies by Theorem \ref{cor-ch}, that $\cO_X(1)$ is $f$-positive. We can apply Theorem \ref{thm: piccolo} to deduce that (1) and (2) hold.
\end{proof}


\subsection{The slope inequality}\label{ssec-slope}
We shall now address the problem of the slope inequality (\ref{slopeinequality}), i.e. of $f$-positivity of the relative canonical sheaf  $\omega_f$.
Let us first establish some preparatory results.

Firstly, we calculate the numerical class of $K_f$ in our setting.
Recall that \cite{miyaoka} the relative canonical bundle of $\pi\colon\pr\to B$ is
\[K_\pi= -rH+\pi^*\det(\mathcal E)\equiv -rH+d\Sigma.\]
Hence, by adjunction theorem, we have that
\begin{equation}\label{eq: canonici}
K_f=\left(K_\pi+\sum X_i\right)_{|X}\equiv \left(\sum k_i-r\right)H_X-\left(\sum y_i-d\right)F=(k_J-r)H_X-(y_J-d)F,
\end{equation}
following Notation \ref{notation-multi} setting $J=\{1,\ldots, c\}$, the whole set of indexes.

Now we recall that $f$-positivity is stable under sum of pullbacks of divisors on the base curve \cite[Sec.2,~Remark 1]{BS3}. It is worth recalling here the proof.
\begin{lemma}
A divisor $D$ on $X$ is $f$-positive if and only if the divisor  $D+f^*M$  is $f$-positive for any divisor $M$ on $B$.
\end{lemma}
\begin{proof}
The $f$-positivity of $D$ is
\[h^0(F,D_{|F}) D^n-n(D_{|F})^{n-1}\deg f_*\mathcal O_X(D)\geq 0
\]
The $f$-positivity of $D+f^*M$ is
\begin{equation}\label{eq: twist}
h^0(F,(D+f^*M)_{|F}) (D+f^*M)^n-n((D+f^*M)_{|F})^{n-1}\deg f_*\mathcal O_X(D+f^*M)\geq 0
\end{equation}
Let us analyze this last inequality. We have:
\begin{itemize}
\item $(D+f^*M)_{|F}\cong D_{|F}$;
\item $(D+f^*M)^n=D^n+nD_{|F}^{n-1}(\deg M)$;
\item $\deg f_*\mathcal O_X(D+f^*M) = \deg f_*\mathcal O_X(D) + (\deg M) \rk f_*\mathcal O_X(D)= \\
\phantom{\deg f_*\mathcal O_X(D+f^*M)}=\deg f_*\mathcal O_X(D) + (\deg M) h^0(F,D_{|F}).$
\end{itemize}
Hence, we have that  \eqref{eq: twist} is indeed
\begin{gather*}
0\leq h^0(F,D_{|F}) (D^n+nD_{|F}^{n-1}(\deg M))-n(D_{|F})^{n-1}\left[\deg f_*\mathcal O_X(D)+(\deg M) h^0(F,D_{|F})\right]=\\ =h^0(F,D_{|F}) D^n-n(D_{|F})^{n-1}\deg f_*\mathcal O_X(D),
\end{gather*}
so, precisely the $f$-positivity of $D$.
\end{proof}
\begin{corollary}\label{cor: equi}
The slope inequality \eqref{slopeinequality} for a morphism $f\colon X\to B$ that  satisfies Assumptions \ref{ass: proper}  is equivalent to the $f$-positivity of $\cO_X(k_J-r)$. If the morphism is balanced then the slope inequality is equivalent to $f$-positivity of $\cO_X(ck-r)$.
\end{corollary}

\begin{theorem}\label{thm: slope}
With the notations above, suppose that $X$ is balanced, and that $k>1$  and that $r<ck$.

Then the following are equivalent:
\begin{enumerate}
\item $K_f^{r-c}\geq 0$;
\item  $\mu(\mathcal E)\geq \frac{\sum y_i}{\sum k_i}=\frac{\sum y_i}{ck}=\frac{y_J}{ck}$.
\end{enumerate}
Moreover, if condition (1) (or (2)) holds, the slope inequality is equivalent to the following inequality:
 \begin{equation}\label{eq: mah}
 (ck-r)h^0(F,K_F)\geq k(r-c)h^0(F_1,K_{F_1}),
 \end{equation} 
 where $f_1\colon Y_1\to B$ is the morphism induced by $\pi $ on the intersection $Y_1$ of all $X_i$'s except for $X_1$, and $F_1$ is a general fibre of $f_1$.
\end{theorem}
\begin{proof}
As in the proof of Theorem \ref{thm: hgg}, let us fix the notation
\[
\alpha:=c\prod_{i=1}^{c} k_i d-\sum_{i=1}^{c}(\prod_{j\not=i}k_j)y_i r=cdk^c-ry_Jk^{c-1}.
\]
Firstly we prove that  $K_f^{r-c}$ is a strict positive multiple of $\alpha$.
This of course tells us that (1) $\iff$ (2).
From (\ref{eq: canonici}) we get
\begin{equation}
K_f^{r-c}=(k_J-r)^{r-c-1}\left((k_J-r)H_X^{r-c}-(r-c)(y_J-d)H_F^{r-c-1}\right)
\end{equation}
Now using the formulas of Lemma \ref{koszul} and dividing by $(k_J-r)^{r-c-1}$ (which is by assumption strictly greater than $0$),
we get
\[
\frac{K_f^{r-c}}{(k_J-r)^{r-c-1}}=d\prod_{i=1}^{c} k_i (k_J-c)-(k_J-r)(\sum_{i=1}^c\prod_{j\not=i}k_j )y_i-(r-c)y_J\prod_{i=1}^{c} k_i.
\]
Now, using the assumption that $X$ is balanced, we have $\prod_{i=1}^{c} k_i=k^c$, $\prod_{j\not=i}k_j =k^{c-1}$, $k_J=ck$ and
we obtain:
\[
\frac{K_f^{r-c}}{(ck-r)^{r-c-1}}=dk^c(ck-c)-(ck-r)(k^{c-1}y_J)-(r-c)k^cy_J =k^{c-1}(k-1)(cdk-ry_J)=(k-1)\alpha.
\]

Let us now turn our attention on the $f$-positivity of $\omega_f$. As noted in Corollary \ref{cor: equi}, we have that this is equivalent to the $f$-positivity of $\mathcal O_X(ck-r)$.

Let us first establish  a formula for the $f$-positivity of $\cO_X(h)$ for any $h\geq 1$ in the balanced case.
\begin{lemma}\label{lem: conti}
Let $h$ be any integer greater or equal to $1$.
Under the assumptions of Theorem \ref{thm: slope}, $f$-positivity of $\cO_X(h)$ is equivalent to:
\begin{equation}\label{eq: contazzo}
\frac{\alpha}{r}\left(h\sum_{i=0}^c(-1)^i \binom{h-ik+r-1}{r-1}\binom{c}{i}+k(r-c)\sum_{i=1}^c(-1)^i \binom{h-ik+r-1}{r-1}\binom{c-1}{i-1}\right)\geq 0.
\end{equation}
\end{lemma}
\begin{proof}\textit{of the Lemma \ref{lem: conti}}.
Let us start from equation \eqref{eq: ine}; $f$-positivity of $\cO_X(h)$ is equivalent to non-negativity of:
\begin{equation}\label{eq: contino}
\alpha\frac{h}{r}\left(\sum_{i=0}^c(-1)^i\sum_{|I|=i}\binom{h-ik+r-1}{r-1}\right)+(r-c)k^c\left(\sum_{i=0}^c(-1)^i\sum_{|I|=i}\binom{h-ik+r-1}{r-1}\frac{ikd-y_I r}{r}\right).
\end{equation}
Let us now observe that the summand on the right of \eqref{eq: contino} can be taken starting from $i=1$, because for $i=0$ we have that $\frac{ikd-y_I r}{r}=0$.
Observe moreover that the left hand side can be simplified via the following equality
\[
\sum_{|I|=i}\binom{h-ik+r-1}{r-1}=\binom{c}{i}\binom{h-ik+r-1}{r-1}.
\]
We now would like to rearrange the term on the right side in equation \eqref{eq: contino} in a similar manner to what we have done in Lemma \ref{lem: asi}.
The key observation is that:
\[
\sum_{i=1}^c(-1)^i\sum_{|I|=i}\binom{h-ik+r-1}{r-1}y_I = \left(\sum_{i=1}^{c}(-1)^i\binom{h-ik+r-1}{r-1}\binom{c-1}{i-1}\right)y_J,
\]
where it is worth recalling that $y_J=\sum_{i=1}^cy_i$.
Using these equalities, we have that   the right side in \eqref{eq: contino} becomes:
\begin{align*}
\sum_{i=1}^c(-1)^i\sum_{|I|=i}\binom{h-ik+r-1}{r-1}\frac{ikd-y_I r}{r}=&\sum_{i=1}^c\frac{(-1)^i}{r}\binom{h-ik+r-1}{r-1}\left(\binom{c}{i}ikd-\binom{c-1}{i-1}y_Jr\right)\\
=&\sum_{i=1}^c(-1)^i\binom{h-ik+r-1}{r-1}\binom{c-1}{i-1}\frac{ckd-y_Jr}{r}\\
=& \frac{\alpha}{rk^{c-1}}\sum_{i=1}^c(-1)^i\binom{h-ik+r-1}{r-1}\binom{c-1}{i-1},
\end{align*}
and the proof is thus concluded.
 \end{proof}
Let us now complete the proof of Theorem \ref{thm: slope} and set  $h=ck-r$.
The expression multiplying $\alpha/r$ in equation \eqref{eq: contazzo} is
\[
(ck-r)\sum_{i=0}^c(-1)^i \binom{ck-r-ik+r-1}{r-1}\binom{c}{i}+k(r-c)\sum_{i=1}^c(-1)^i\binom{ck-r-ik+r-1}{r-1}\binom{c-1}{i-1}.
\]
Now observe that 
\[
\sum_{i=1}^c(-1)^i \binom{ck-r-ik+r-1}{r-1}\binom{c-1}{i-1}=-\sum_{j=0}^{c-1}(-1)^{j}\binom{(c-1)k-r-jk+r-1}{r-1}\binom{c-1}{j}.
\]
Now, let us call $Y_1$ be the intersection in $\mathbb P$ of $X_2,\ldots, X_c$, and let $f_1\colon Y_1\to B$ be the morphism induced.
We see from the very computation made in Proposition \ref{koszul} that the above expression is
$-\rk (f_1)_*\cO_{Y_1}((c-1)k-r)$.
So, the slope inequality is equivalent to the non-negativity of
\begin{align*}
\frac{\alpha}{r}\left[(ck-r)\rk f_*\cO_X(ck-r)-k(r-c)\rk (f_1)_*\cO_{Y_1}((c-1)k-r)\right]=\\
=\frac{\alpha}{r}\left[(ck-r)h^0(F,K_F)-  k(r-c)h^0(F_1,K_{F_1})\right].
\end{align*}
\end{proof}
\begin{remark}
Note that condition $r<ck$ in Theorem \ref{thm: slope} implies that the canonical sheaf on the fibres of $f$ is very ample, i.e. that $K_f$ is relatively very ample, so the fibres are of general type.
\end{remark}
\begin{remark}
In the case $c=1$, we have that $F_1=\mathbb P^{r-1}$ and so inequality \eqref{eq: mah} is trivially satisfied, so   \cite[Theorem 1.2 (1)]{B-S-hyper} is implied by this result.
\end{remark}
\begin{remark}
In the non balanced case, one could still obtain an inequality, but more involved.
\end{remark}
In general it does not seem clear to understand the positivity of  expression \eqref{eq: mah}. 
However, we can prove the slope inequality in some cases  as follows:
\begin{proposition}\label{prop: meglio che niente}
Under the assumptions of Theorem \ref{thm: slope}, suppose that  one of the following conditions holds:
\begin{enumerate}
\item $(c-1)k<r$ (this is equivalent to asking that $F_1 $ is not of general type). 
\item  $c=2,3,4$  and $k\gg0$.
\item $c$ is fixed and $r\gg0$.
\end{enumerate}
Then the slope inequality is equivalent to  $\mu(\mathcal E)\geq \frac{y_J}{ck}$.
\end{proposition}
\begin{proof}
In case (1) $h^0(F_1,K_{F_1})$  is trivially zero, so by \eqref{eq: mah}  the slope inequality is equivalent to the non-negativity of $\alpha$, as wanted.

As for case (2), let us compute \eqref{eq: mah} as a polynomial in $k$, and see its leading coefficient.
\begin{align*}
(ck-r)\sum_{i=0}^c(-1)^i \binom{ck-r-ik+r-1}{r-1}\binom{c}{i}+k(r-c)\sum_{i=1}^c(-1)^i\binom{ck-r-ik+r-1}{r-1}\binom{c-1}{i-1}\\
=(ck-r)\binom{ck-r+r-1}{r-1}+\sum_{i=1}^c (-1)^i\binom{ck-r-ik+r-1}{r-1}\left( (ck-r)\binom{c}{i}+k(r-c)\binom{c-1}{i-1}\right)\\
=\frac{k^r}{(r-1)!}\left[ c^r+\sum_{i=1}^c(-1)^i(c-i)^{r-1}\left(\frac{c^2}{i}+r \right)\binom{c-1}{i-1}\right]+\cO(k^{r-1}).
\end{align*}
Now, if we compute the coefficient for small $c$, we can see that it is non-negative:
\begin{itemize}
\item For $c=2$ we have $2^r-(r+4)$, which is always strictly greater than $0$ as soon as $r\geq 2$ (indeed in our case $r$ need to be greater than $4$).
\item For $c=3$ we obtain $3^r-2^{r-1}(r+9)+2(r+9/4)$, which is strictly greater than $0$  for $r\geq 5$.
\item For $c=4$ we obtain $4^r-3^{r-1}(r+16)+2^{r-1}3(r+4)-3(r+16/9)$, which is greater than $0$ for $r\geq 6$
\end{itemize}

As for the third case, we just need to observe that the leading coefficient as a polynomial in $k^r$ in the above expression is positive if $c$ is fixed and $r\gg0$.
\end{proof}

\begin{remark}
Note that the result (2) for $k\gg 0$ can not be derived from $f$-positivity of $\cO_X(h)$ for $h\gg 0$, because changing $k$ does change also $X$.
\end{remark}

\begin{remark}\label{rem: confronto}
In \cite[Lemma 1.1]{ENO} Enokizono proved the following: in the balanced case with $c=r-2$ (so $X$ is a surface), the following equalities hold:
\begin{align*}
&K_f^2=((k-1)r-2k)(k-1)\alpha;\\
&\deg f_*\omega_f = \frac{1}{24}((3r-5)k-3r+1))(k-1)\alpha.
\end{align*}
So, we see that the non-negativity of both these invariants correspond to the non-negativity of $\alpha$, coherently with our results.
It is easy to see that this equalities implies the slope inequality in case $\alpha\geq 0$. In particular, by Enokizono's result, there is always an equality, regardless of the sign of $\alpha$:
\[
K_f^2=\frac{24((r-2)k-r)}{(3r-5)k-3r+1}\deg f_*\omega_f.
\]
\end{remark}


\section{The cones of cycles of $\pr_B(\cE)$ and an instability condition}\label{birational}
\subsection{The cone of ``$f$-positive complete intersections'' in $N^c(\pr)$}

In this section we interpret some of the results obtained in terms of cones in the real N\'eron-Severi space of codimension $c$ cycles of $X$.
We see that inequality \eqref{eq: principale} defines a meaningful cone, which is intermediate with respect to the nef and the pseudoeffective ones.
For stating this, we rewrite in a more compact form a result of Fulger \cite{fulger}, that describes completely this last two cones.

Recall first that given any vector bundle $\cE$ over a curve $B$, the real N\'eron-Severi space of codimension $c$ cycles of $\pr=\pr_B(\cE)$
\[N^c(\pr) := \frac{\{\mbox{real span of classes of $c$-dimensional subvarieties of $\pr$}\}}{\mbox{numerical equivalence}}\]
is 2-dimensional, generated by the classes $H^c$ and $H^{c-1}\Sigma$, where $H=[\cO_\pr(1)]$ and $\Sigma$ is the class of a fibre of $\pi$.

Let us consider the Harder-Narasimhan filtration of $\cE$
\[0=\cE_0\subset \cE_1\subset \ldots \subset \cE_l=\cE,\]
and call $\mu_i:= \mu(\cE_i/\cE_{i-1})$, and $\mu:=\mu(\cE)$. Recall that $r_i:= \rk \cE_i$, and that $\rk\cE=r=r_l$.
Recall in particular that
\begin{equation}\label{HN}
\mu_{\ell}< \mu_{\ell-1} <\ldots< \mu_1.
\end{equation}

It is useful to introduce the following notation.
\begin{definition}
With the above notations, we define
the {\em virtual slopes} ${\widetilde\mu}_1\geq{\widetilde\mu}_2\geq \ldots \geq {\widetilde\mu}_r$ of $\cE$ as follows.
Let $r_j$ be the rank of $\cE_j$ and let $i\in \{1,...,r-1\}$. If $r_j \leq i < r_{j+1}$ then we define ${\widetilde \mu}_i=\mu_j$. We define coherently ${\widetilde \mu}_r=\mu_{\ell}$.
\end{definition}
 Observe that $d={\rm deg}\cE=r\mu=\sum_{i=1}^r{\widetilde \mu}_i$.
With this notation, we can restate Fulger's result as follows.

\begin{theorem}[Fulger \cite{fulger}]
With the notations defined above, we have:
\begin{equation}\label{psefc}
{\rm Pseff}^c(\pr)= \left< H^{c-1}\Sigma\,,\,H^{c}-(\widetilde\mu_1+\ldots +\widetilde\mu_c)  H^{c-1}\Sigma\right>;
\end{equation}
\begin{equation}\label{nefc}
{\rm Nef}^c(\pr)= \left< H^{c-1}\Sigma\,, \,H^{c}-(\widetilde\mu_{r-c+1}+\ldots +\widetilde\mu_r)  H^{c-1}\Sigma\right>.
\end{equation}
\end{theorem}
\begin{proof} The first formula is just the content of Theorem 1.1 in \cite{fulger} with the notation suitably adapted.
Observe that in \cite{fulger} the indexes of slopes are reversed with respect to our definition.
With this in mind, and following the construction given by Figure 2 in \cite{fulger} and the notation therein, one immediately obtains that $\nu ^{(i)}=-(\widetilde\mu_1+\ldots +\widetilde\mu_c)$.

For the second formula one just needs to use the fact that  ${\rm Nef}^c(\pr)$ is the subset of ${\rm Pseff}^c(\pr)$ defined by positivity product with
\[{\rm Pseff}_c(\pr)={\rm Pseff}^{r-c}(\pr)=<H^{r-c-1},H^{r-c}-({\widetilde \mu}_1+\ldots+{\widetilde \mu}_{r-c})\Sigma>,\]
 and so it is the two dimensional cone determined by $H^{c-1}\Sigma$ and $H^c-aH^{c-1}\Sigma$ such that
\[0=\left( H^c-aH^{c-1}\Sigma \right)\left(H^{r-c}-({\widetilde \mu}_1+\ldots+{\widetilde \mu}_{r-c})\Sigma\right)=H^r-(a+{\widetilde \mu}_1+\ldots+{\widetilde \mu}_{r-c})=\]
\[=d-(a+{\widetilde \mu}_1+\ldots+{\widetilde \mu}_{r-c})=\widetilde\mu_{r-c+1}+\ldots +\widetilde\mu_r-a.\]

\end{proof}
We can now reformulate the results of the previous section using the language of cones.
\begin{definition}\label{def: cono}
Let $\mathbb B$ be the cone in ${\rm Nef}^c(\pr)$ generated by the classes $[H^{c-1}\Sigma]$ and $[H^{c}-c\mu  H^{c-1}\Sigma]$.
\end{definition}
Note that for any $c\in \{1,\ldots,r-1\}$ we have that
\begin{itemize}
\item $\sum_{i=1}^c\widetilde \mu_i>c\mu$;
\item $\sum_{i=1}^c\widetilde \mu_{r-i+1} < c\mu$.
\end{itemize}
This means that the cone $\mathbb B$ is indeed contained in the pseudoeffective cone and contains the nef cone:
\[
{\rm Nef}^c(\pr)\subseteq \mathbb B\subseteq {\rm Pseff}^c(\pr).
\]
Note that by Theorem \ref{thm: M-N}  the inclusions are strict unless $\cE$ is semistable (in which case the cones all coincide).
Theorems \ref{thm: piccolo}, \ref{thm: hgg} and \ref{thm: slope} tell us the following:
\begin{proposition}\label{cone}
In the above notations, let $X\subset \pr$ be a codimension $c$ cycle which is a complete intersection of $c$ relative hypersurfaces $X_1,\ldots, X_c$ in $\pr$ of degree at least 2.
\begin{enumerate}
\item The numerical class of $X$ is contained in $\mathbb B$ if and only if there exists $h< \min\{k_i\}$ such that  $\cO_X(h)$ is $f$-positive;
\item the numerical class of $X$ is contained in $\mathbb B$ if and only if  $\cO_X(h)$ is $f$-positive for any $h< \min\{k_i\}$;
\item If $X$ is balanced and  the numerical class of $X$ lies in the interior of $\mathbb B$, then  $\cO_X(h)$ is strictly $f$-positive for $h\gg 0$;
\item If $X$ is balanced and  $\cO_X(h)$ is $f$-positive for $h\gg 0$, then the numerical class of $X$ is contained in $\mathbb B$;
\item If $X$ is balanced, $r<kc$ and (1) or (2) in Proposition \ref{prop: meglio che niente} hold, then the numerical class of $X$ is contained in $\mathbb B$ if and only if the slope inequality holds.
\end{enumerate}
\end{proposition}


\subsection{An instability condition for of the fibres}
We can use the results of Section \ref{sec-inequalities} vice versa  as in \cite{B-S-hyper}, and prove the following {\em instability} condition for the fibres of a global relative complete intersection. As usual, let $\cE$ be a rank $r\geq 3$ vector bundle over a curve $B$, and let $\pi\colon\pr=\pr_B(\cE)\to B$ be the projective bundle.

\begin{corollary}\label{instability}
Let $X\subset \pr$ a relative complete intersection in the projective bundle $\pr=\pr_B(\cE)$ satisfying Assumptions \ref{ass: proper}.
If $\sum_{i=1}^c\frac{y_i}{k_i}>c\mu$ (equivalently $[X]\not\in \mathbb B$), then:
\begin{itemize}
\item[(i)] the fibres of $f$ are Chow unstable with the restriction of $\cO_{\mathbb P ^{r-1}}(h)$ for any $h< \min\{k_i\}$.
\item[(ii)] Assume moreover that $X$ is balanced. Then the fibres of $f$ are Chow unstable with the restriction of $\cO_{\mathbb P ^{r-1}}(h)$ for  $h\gg0$;
\item[(iii)] Assume moreover that $X$ is balanced,  $r<kc$ and (1), (2) or (3) in Proposition \ref{prop: meglio che niente} holds. Then the fibres of $f$ are unstable with respect to their dualizing sheaf.
\end{itemize}
\end{corollary}
\begin{proof}
Immediate from Theorem \ref{thm: main} and Theorem \ref{teo-ch}.
\end{proof}

\begin{remark}
In \cite{B-S-hyper}, in the codimension one case, we proved a more general instability condition, and this led us to a singularity condition (a bound on the log canonical threshold of the fibres of $f$ via Lee's result).
In the general codimension case, of course from  Lee's result we could obtain a singularity condition on the Chow form of the fibres; but as for the fibres themselves, it is not so easy to get geometric information from an instability condition.
\end{remark}
As an application of these last results, together with the detailed study of the hypersurface case made in \cite{B-S-hyper}, we 
use  Corollary \ref{instability} to construct families of complete intersections whose general fibre is of general type, asymptotically unstable and has only one (very) singular point, as follows. 
\begin{proposition}\label{esempio}
Let $a$ be a positive natural number, and let $r\geq 3$ be a natural number. Consider the rank $r\geq 3$ vector bundle over $\pr^1$
\[\cE:=\cO_{\pr^1}(a)^{\oplus r-1}\oplus \cO_{\pr^1}(a-1).\]
Let $0<c<r-1$ be an integer.
Then there exists an $m\in \mathbb N$ big enough, such that calling  $X$ the complete intersection of $c$ general members of the linear system $$X_i\in |m((r+1)
a-1)H-m(r+1)\Sigma|$$ on
$\pr_B(\cE)$ and $f\colon X\to B$ the induced morphism, we have that:
\begin{enumerate}
\item the general fibre $F$ of $f$ is a complete intersection of general type with only one singular point;
\item the couple $(F, hH_{|F})$ is Chow unstable for any $h\in \mathbb N^{>0}$;
\end{enumerate}
\end{proposition}
\begin{proof}
The Harder-Narasimhan sequence of $\cE$ is simply
\[
\xymatrix{
\cE_0\ar@{^{(}->}[r]&\cE_1 \ar@{^{(}->}[r]&\cE_2\\
0\, \ar@{=}[u]\ar@{^{(}->}[r]&\,\,\,\cO_{\pr^1}(a)^{\oplus r-1}\,\ar@{=}[u] \ar@{^{(}->}[r]&\cE\ar@{=}[u]}
\]
and $\mu_1=a>\mu(\cE)=\frac{(r-1)a+a-1}{r}=a-\frac{1}{r}>\mu_2=a-1$.
We have that 
\[
\frac{m[((r+1)a-1)]}{m(r+1)}=\frac{(r+1)a-1}{r+1}=a-\frac{1}{r+1}, 
\] so, as this ratio is between $\mu_1$ and $\mu_2$:
\[\mu_1=a>a-\frac{1}{r+1}>a-1=\mu_2,\] we are in the conditions of applying \cite[Theorem 1.7]{B-S-hyper}. 
That theorem states that for $m\gg 0$ the general member of  $|m((r+1)a-1)H-m(r+1)\Sigma|$ has $\pr(\cE/\cE_1)$ as base locus, so if we consider $c$ general members of this linear system, they will intersect in a variety $X$ smooth outside the section $\pr(\cE/\cE_1)$. Now we just need to observe that
\[a-\frac{1}{r+1}>\mu(\cE),\] so by Corollary \ref{instability} we obtain the statement.
\end{proof}
\begin{remark}\label{rem-referee}
Note that an asymptotically Chow unstable complete intersection of general type need to be singular. Indeed, let $X$ be such a variety. By \cite{mabuki}, Chow-stability for $h\gg 0$ is equivalent to Hilbert-stability for $h\gg 0$. By \cite[Corollary 4]{don}, this is implied by the existence of a K\"ahler-Einstein metric on $X$ (because the automorphism group of $X$ is finite). Hence, if a smooth complete intersection is not Chow-stable for $h\gg 0$, then it cannot carry a K\"ahler-Einstein metric.
When $K_X$ is ample this cannot happen by Aubin-Yau Theorem.
\end{remark}

\subsection{Acknoledgements} 
We wish to thank Yongnam Lee for enlightening conversations, and in particular for having pointed out to us the results of Ferretti used for Proposition \ref{stab-intersection}.
Some anonymous referees helped to improve a lot this paper, by pointing out several inaccuracies and some mistakes, that have now been fixed. Moreover, we literally owe to one of them Remark \ref{rem-referee}. We thank them heartily.


\medskip

\noindent {Miguel \'Angel Barja \\ Departament de Matem\`atiques \\ Universitat Polit\`ecnica de Catalunya-BarcelonaTech
\\ Institut de Matemàtiques de la UPC-BarcelonaTech (IMTech) \\Centre de Recerca Matemàtica (CRM)
\\Avda. Diagonal 647, 08028 Barcelona, Spain.}\\
E-mail: \textsl {miguel.angel.barja@upc.edu}

\medskip
\noindent Lidia Stoppino,\\Dipartimento di Matematica, Universit\`a di Pavia,\\ Via Ferrata 5, 27100 Pavia, Italy.\\
E-mail: \textsl {lidia.stoppino@unipv.it}.


\begin{thebibliography}{99}

\bibitem{ACG2} E.~Arbarello, M.D.T.~Cornalba, P.A.~Griffiths, {\em Geometry of algebraic curves. Vol II}. Grundlehren der mathematischen Wissenschaften, 268. Springer, Heidelberg, 2011. xxx+963 pp.


\bibitem{barja-preprint} M.A.~Barja, {\em Higher Dimensional Slope Inequalities for Irregular fibrations}, to appear in Annali della Scuola Normale di Pisa.
\textit{DOI 10.2422/2036-2145.202109-012}

\bibitem{Barja-BUMI} M.A.~Barja, {\em Slope inequalities for fibrations of non-maximal Albanese dimension.} Boll. Unione Mat. Ital. 15 (2022), no. 1-2, 315.

\bibitem{Barja3folds} M.A.~Barja, \emph{On the slope of fibred threefolds}, Internat. J. Math. {\bf 11} n.4 (2000), 461-491.

\bibitem{BS3} M.A.~Barja, L.~Stoppino, {\em Stability conditions and positivity of invariants of fibrations}, Algebraic and Complex Geometry,
Springer Proc. in Math. \& Stat. Vol. 71, 2014, 1--40, volume in honour of Klaus Hulek's 60th birthday.

\bibitem{B-S-hyper} M.A.~Barja, L.~Stoppino, {\em Stability and singularities of relative hypersurfaces}, IMRN, no. 4 (2016), 1026--1053.

\bibitem{bost} J.~Bost, {\em Semi-stability and heights of cycles}, Invent. Math.,{\bf 118 (2)}  (1994), 223-253.

\bibitem{CTV} G.~Codogni, L.~Tasin, F.~Viviani, {\em Slope inequalities for KSB-stable and K-stable families},  Proc. Lond. Math. Soc. (3) 126 (2023), no. 4, 1394--1465. 14D06 (14J10)

\bibitem{C-H} M.~Cornalba, J.~Harris, \emph{Divisor classes associated to families of stable varieties, with applications to the moduli space of curves.} Ann. Sc. Ec. Norm. Sup. {\bf 21 (4)} (1988), 455-475.

\bibitem{don} S.K.~Donaldson, {\em Scalar curvature and projective embeddings. I.} J. Differential Geom. {\bf 59} (2001), no. 3, 479--522. 

\bibitem{E-H} D.~Eisenbud, J.~Harris, \emph{The geometry of schemes}. Graduate Texts in Mathematics, 197. Springer-Verlag, New York, 2000. x+294 pp.

\bibitem{ENO} M.~Enokizono, {\em Dufree-type inequality for complete intersection surface singularities}, Duke Math. J. {\bf 170 (1)} (2021), 1--21.

\bibitem{ferretti} R.G.~Ferretti, \emph{Diophantine approximation and toric deformation}, Duke Mat. J. {\bf 118 (3)} (2003), 493--522.

\bibitem{fulger} M.~Fulger, {\em The cones of effective cycles on projective bundles over curves}, Math. Z. {\bf 269} (2011), no. 1-2, 449--459.

\bibitem{fulton} W.~Fulton, { Intersection Theory, second edition}, Springer-Verlag, 1998.

\bibitem{GKM} A. Gibney, S. Keel, I. Morrison, {\em Towards the ample cone of {$\overline M_{g,n}$}}, {J. Amer. Math. Soc.} {\bf 15} (2002), no. {2}, {273--294}.

\bibitem{H-N} G.~Harder, M.S.~Narasimhan, , {\em On the cohomology groups of moduli spaces of vector bundles on curves}, Math. Annalen {\bf 212}  (1975), no. 3, pp.  215--248,

\bibitem{Har} R.~Hartshorne, Algebraic geometry. GTM, No. 52. Springer-Verlag, New York-Heidelberg, 1977.

\bibitem{zhang} Y.~Hu, T.~Zhang, {\em Fibered varieties over curves with low slope and sharp bounds in dimension three.} J. Algebraic Geom. {\bf 30} (2021), no. 1, 57--95. 

\bibitem{H-Z} Y.~Hu, T.~Zhang, {\em Relative Severi inequality for fibrations of maximal Albanese dimension over curves.} Forum Math. Sigma {\bf 10} (2022), Paper No. e45, 31 pp.

\bibitem{KM} J. Kollar and S. Mori, Birational geometry of algebraic varieties. Cambridge Tracts in Mathematics,134. Cambridge University Press, Cambridge, 1998. viii+254 pp.

\bibitem{Lee}  Y. Lee, {\em Chow stability criterion in terms of log canonical threshold}, J. Korean Math. Soc., {\bf 45 (2)} (2008), 467--477.

\bibitem{mabuki} T.~Mabuchi, {\em Chow-stability and Hilbert-stability in Mumford's geometric invariant theory.} Osaka J. Math. 45 (2008), no. 3, 833--846. 

\bibitem{miyaoka} Y.~Miyaoka, {\em The Chern class and Kodaira dimension of a minimal variety}, Algebraic Geometry, Sendai (1985), 449--476.

\bibitem{mum} D.~Mumford, {\em Stability of projective varieties.} Enseign. Math. (2) {\bf 23 }(1977), no. 1-2, 39--110.

\bibitem{nak} N.~Nakayama,  {\em Zariski-decomposition and abundance.} MSJ Memoirs, 14. Mathematical Society of Japan, Tokyo, 2004. xiv+277 pp.


\bibitem{ohno} K.~Ohno, \emph{Some inequalities for minimal fibrations of surfaces of general type over curves}, J. Math. Soc. Japan {\bf 44 (4)} (1992), 643-666.

\bibitem {pardini} R.~Pardini, {\em The Severi inequality $K^2\geq 4\chi$ for surfaces of maximal Albanese dimension}, Invent. Math. {\bf 159} 3 (2005), 669--672.

\bibitem{LS} L.~Stoppino, {\em Slope inequalities for fibered surfaces via GIT},  Osaka J. Math. Vol. 45, No. 4 (2008), 1027--1041.

\bibitem{xiao} G.~Xiao, \emph{Fibred algebraic surfaces with low slope.} Math. Ann. {\bf 276} (1987), 449-466.


\end{thebibliography}
\end{document}